\documentclass[11pt]{article}%{smfart}%
\usepackage{amssymb}\usepackage{amsmath}
\newcommand{\be}{\begin{eqnarray}}
\newcommand{\ee}{\end{eqnarray}}
\newcommand{\nn}{\nonumber}
\newcommand{\tr}{\mathop{\rm tr}\nolimits}
\newcommand{\one}{{\mathfrak 1}}
\newcommand{\two}{{\mathfrak 2}}
\newcommand{\three}{{\mathfrak 3}}

\newcommand{\proof}{{\em Proof.\ }}

\def\qed{\hfill $\square$}
\def\PP{\mathbb P}
\def\EE{\mathbb E}
\def\la{\lambda}
\def\SF{\mathsf F}
\def\SG{\mathsf G}
\def\SH{\mathsf H}
\newcommand{\rf}[1]{(\ref{#1})}

\renewcommand{\=}[1]{\stackrel{(\ref{#1})}{=}}
\newcommand{\FRW}[6]{\Bigl\{\begin{smallmatrix}%
 \textstyle #1\, & \textstyle #3\, & \textstyle #5 \\[2pt]%
 \textstyle #2\, & \textstyle #4\, & \textstyle #6 %
 \end{smallmatrix}\Bigr\}}

\newtheorem{prop}{Proposition}
\newtheorem{lem}{Lemma}
\begin{document}
% June 2006  
\strut\hfill math.QA/0609519 \\
\vspace*{3mm}
\begin{center}
{\large\bf On one ansatz for $sl_2$--invariant R--matrices } \\
[3mm] {\sc A. G. Bytsko} \\ [2mm]
{ Steklov Mathematics Institute \\
 Fontanka 27, 191023, St.-Petersburg, Russia }
\footnotetext{This work was supported in part by the grant INTAS
 YS--03-55-962, the Russian Fund for Fundamental Research grant 05-01-00922,
 and by a grant from the Russian Science Support Foundation.}

\end{center}
\vspace{1mm}
\begin{abstract}

The spectral decomposition of regular $sl_2$--invariant R--matrices $R(\la)$
is studied by means of the method of reduction of the Yang--Baxter equation
onto subspaces of a given spin. Restrictions on the possible structure of
several highest coefficients in the spectral decomposition are derived. The
origin and structure of the exceptional solution in the case of spin $s=3$ 
are explained. Analogous analysis is performed for constant R--matrices. 
In particular, it is shown that the permutation matrix $\PP$ is a ``rigid''
solution.

\end{abstract}
\section*{\S 1. Introduction}   %%%%%%%%%%%%%%%%%%%%%%%%%%%%%%%%%%%

The Yang--Baxter equation plays a key role in the quantum inverse 
scattering method (see, e.g., the reviews \cite{KS,Fad1}). Its
braid group form looks as follows
\be\label{YB}
 R_{\one\two}(\lambda)\, R_{\two\three}(\lambda+\mu) \,
 R_{\one\two}(\mu) =  R_{\two\three}(\mu)\,
 R_{\one\two}(\lambda+\mu)\, R_{\two\three}(\lambda) \,.
\ee 
In this article we will consider the Yang--Baxter equation \rf{YB} on the 
space $V_{s}^{\otimes 3}$, where $V_s$ is an irreducible finite--dimensional
representation of the algebra~$sl_2$. The dimension of the representation
$V_s$ is $(2s{+}1)$, where $s$ is a positive integer or semi--integer
number (referred to below as {\em spin}). Here and below we use the 
standard notations: the lower indices of $R(\lambda)$ indicate the tensor 
components of~$V_{s}^{\otimes 3}$ where $R(\lambda)$ acts nontrivially.

An operator--valued function
$R(\lambda) : {\mathbb C} \mapsto {\rm End}\ V_{s}^{\otimes 2}$
that satisfies \rf{YB} is called an \hbox{R--matrix}. We will consider
$sl_2$--invariant R--matrices, i.e., those that have the spectral
decomposition of the form
\be\label{spec}
 R(\la) = \sum_{j=0}^{2s}  r_j(\lambda) P^{j} \,.
\ee
Here $P^j$ is the projector onto $V_j$ which is the subspace of
spin $j$ in~$V_s^{\otimes 2}$, and $r_j(\la)$ is a scalar function.
Additionally, we assume that R--matrices under consideration
are  {\em regular}, {\em unitary}, and normalized, that is, the
following relations are satisfied
\be\label{run}
 r_j(0) = 1 \,, \qquad r_j(\la) r_j(-\la) =1 \,,
 \qquad r_{2s}(\la) = 1\,.
\ee
Let us remark that unitarity is a consequence of
regularity and normalization~\cite{By}.

Since regular R--matrices can be used to construct {\em local}
integrals of motion for lattice models, in particular for
spin chains, the problem of finding all solutions of the
Yang--Baxter equation satisfying the properties \rf{run} is 
important for the quantum inverse scattering method. At present,
there are known four series of inequivalent $sl_2$--invariant 
regular solutions and one exceptional solution for $s=3$ 
(see~\cite{Ke} and references therein). A computer--based 
check~\cite{Ke} led to a conjecture that this list of 
solutions is exhausting. However, the corresponding
classification theorem has not been proven yet.
In the present article, applying the approach developed in \cite{By},
we will make some progress in this direction. In particular,
we will explain the origin and the structure of the 
exceptional solution for $s=3$.

The paper is organized as follows. \S2 contains analysis of
one ansatz for an R--matrix. Although the results presented
here are well known, we provide all necessary technical details
because our aim is to develop similar technique in a more 
general case. In~\S3 we remind briefly the main details
of the approach developed in \cite{By} for analysis of $sl_2$-- 
and $U_q(sl_2)$--invariant R--matrices. Here we also prove
one useful additional relation (Lemma~\ref{PM}). In~\S4.1 we will 
demonstrate that analysis of some number of highest coefficients in 
the spectral decomposition of an R--matrix can be done in a way
closely resembling the analysis described in~\S2.
In~\S\S4.2--4.3 we will give details of this analysis.
In particular, it turns out that the exceptional solution 
arises as a consequence of degeneration of a certain set
of matrices. In~\S5 we perform analogous analysis for constant
R--matrices. In particular, it is shown that the permutation
$\PP$ is a ``rigid'' solution. The Conclusion summarizes the
main results.

\section*{\S 2. Analysis of one ansatz for an R--matrix} %%%%%%%%%%%%%%%%

Let $\mathbb E$ denote the identity operator on~$V_{s}^{\otimes 2}$.
For $s\geq 1$ let us consider R--matrices of the following form
\be\label{ans}
 R(\lambda) = \frac{1}{1+ f(\la)} \, \bigl( {\mathbb E} +
 f(\la) \, {\mathbb P} + g(\lambda) \, P^0 \bigr) \,.
\ee
Here $\PP$ is the permutation operator on $V_{s} \otimes V_{s}$.
Recall that it can be expressed in terms of projectors:
\be\label{perm}
 \PP = \sum_{j=0}^{2s} (-1)^{2s-j} \, P^j  \,.
\ee
If the scalar functions $f(\la)$ and $g(\la)$ satisfy the
condition $f(0)=g(0)=0$, then \rf{ans} is an ansatz for 
a solution of the Yang--Baxter equation in the class~\rf{run}. It turns
out that all R--matrices of this type can be described explicitly.

\begin{lem}\label{PPrel1}
The following relations hold on~$V_{s}^{\otimes 3}$
\be
\label{pp1}
 P^0_{l}  P^0_{l} = P^0_{l} \,,\quad \PP_{l} \, \PP_{l} = \EE \,, &&
 P^0_{l} \, \PP_{l} = \PP_{l} \, P^0_{l} = \xi \, P^0_{l} \,, \\
\label{pp2}
  \PP_{l} \, \PP_{l'} \, \PP_{l} &=& \PP_{l'} \, \PP_{l} \, \PP_{l'} \,, \\
\label{pp3}
 P^0_{l} \, \PP_{l'} \, \PP_{l} = \PP_{l'} \, \PP_{l} \, P^0_{l'} \,, &&
 \PP_{l} \, P^0_{l'} \, \PP_{l} = \PP_{l'} \, P^0_{l} \, \PP_{l'} \,, \\
\label{pp4}
 P^0_{l} \, \PP_{l'} \, P^0_{l} = \eta \, P^0_{l} \,, &&
 P^0_{l} \, P^0_{l'} \, P^0_{l} = \eta^2 \, P^0_{l} \,, \\
\label{pp5}
 P^0_{l} \, P^0_{l'} \, \PP_{l} = \xi \eta \, P^0_{l} \, \PP_{l'} \,, &&
 \PP_{l} \, P^0_{l'} \, P^0_{l} = \xi \eta \, \PP_{l'} \, P^0_{l} \,,
\ee
where $l={\scriptstyle \{12\}}$, $l'={\scriptstyle \{23\}}$ or
$l={\scriptstyle \{23\}}$, $l'={\scriptstyle \{12\}}$, and
$\xi$ and $\eta$ are scalar constants:
\be\label{xieta}
 \xi = (-1)^{2s} \,, \qquad \eta = \frac{1}{2s+1} \,.
\ee
\end{lem}

\proof
The third relation in \rf{pp1} follows from~\rf{perm}. Equalities
\rf{pp2} and \rf{pp3} are obvious. Relations \rf{pp4} follow from
the well--known relation (see, e.g.~\cite{By})
\be\label{ppp}
  P^0_{\one\two} \, P^j_{\two\three} \, P^0_{\one\two} =
    \frac{2j+1}{(2s+1)^2} \, P^0_{\one\two}   \,.
\ee
Relation \rf{pp5} can be derived as follows:
\be\label{pp5pr}
 P^0_{\one\two} \, P^0_{\two\three} \, \PP_{\one\two} =
 P^0_{\one\two} \,  \PP_{\one\two} \, P^0_{\one\three}
 \={pp1} \xi\, P^0_{\one\two} \, P^0_{\one\three} =
 \xi\, P^0_{\one\two} \, \PP_{\two\three} \,
    P^0_{\one\two} \, \PP_{\two\three} \={pp4}
 \xi \eta \, P^0_{\one\two} \, \PP_{\two\three} \,.
\ee
\qed\\[-0.5mm]

Substituting \rf{ans} in \rf{YB} and using the relations of 
Lemma~\ref{PPrel1}, it is not difficult to check that the 
Yang--Baxter equation for the ansatz under consideration is equivalent
to the following equation
\begin{equation}\label{FGH}
 F_{\la,\mu} \, \SF +
 G_{\la,\mu} \, \SG  +
 H_{\la,\mu} \, \SH  +
 H_{\mu,\la} \, \tilde{\SH} = 0 \,,
\end{equation}
where
\begin{equation}\label{FGH0}
\begin{aligned}
\SF &= \PP_{\one\two} - \PP_{\two\three} \,, \quad &
\SG &= P^0_{\one\two} - P^0_{\two\three} \,, \\
\SH &= P^0_{\one\two} \, \PP_{\two\three} -
    \PP_{\one\two} \, P^0_{\two\three} \,,
 \quad & \tilde{\SH} &= \PP_{\two\three}\, P^0_{\one\two} -
    P^0_{\two\three}\, \PP_{\one\two}  \,,
\end{aligned}
\end{equation}
and
\be
\label{F}
 F_{\la,\mu} &=& f(\la) + f(\mu) - f(\la+\mu) \,, \\
\label{G}
 G_{\la,\mu} &=& g(\la) + g(\mu) - g(\la+\mu) + \xi f(\la) g(\mu)
 + \xi g(\la) f(\mu) \\
\nn && +\ g(\la) g(\mu) + \eta g(\la) g(\mu) f(\la+\mu)
 + \eta^2 g(\la) g(\mu) g(\la+\mu) \,, \\
\label{H}
 H_{\la,\mu} &=& g(\la) f(\la+\mu) - f(\la) g(\la+\mu)
 + \xi \eta \, g(\la) f(\mu) g(\la+\mu) \,.
\ee

\begin{lem}\label{LD}
 For $s\geq 1$ the matrices $\SF$, $\SG$, $\SH$ and $\tilde{\SH}$
 in~\rf{FGH0} are linearly independent.
\end{lem}
The proof is given in Appendix~B.

As a consequence of Lemma~\ref{LD} it follows that equations \rf{FGH} 
are equivalent to the following system of functional equations
\be \label{F1}
 F_{\la,\mu} &=& 0 \,, \\
\label{G1}
 G_{\la,\mu} &=& 0 \,, \\
\label{H1}
 H_{\la,\mu} = H_{\mu,\la} &=& 0 \,.
\ee

Analysis of system \rf{F1}--\rf{H1} is fairly simple.
There are three non--trivial cases:

1) $f(\la)\neq 0$, $g(\la)=0$. In this case it is obvious from
\rf{F} that $f(\la)$ is a linear function. Without loss of
generality one can choose $f(\la)=\la$.

2) $f(\la)=0$, $g(\la) \neq 0$. In this case there remains
one equation on $g(\la)$,
\be\label{G2}
 g(\la) + g(\mu) - g(\la+\mu)  + g(\la) g(\mu)
 + \eta^2 g(\la) g(\mu) g(\la+\mu) = 0 \,,
\ee
which has (for $\eta \neq \frac{1}{2}$) the following solution:
\be\label{gsol1}
 g(\la)=
  b\, \frac{1-e^{\gamma\la}}%
  {e^{\gamma\la} - b^2} \,, \qquad
 \quad b+b^{-1}=\eta^{-1} \,.
\ee
Here $\gamma$ is an arbitrary finite constant, which can be
chosen unity without loss of generality. 

3) $f(\la) \neq 0$, $g(\la) \neq 0$. One can again choose 
$f(\la)=\la$. Introducing a new function $h(\la)=f(\la)/g(\la)$,
we can rewrite equations \rf{H1} in the following form
\be\label{H2}
 h(\la+\mu) = h(\la) - \xi \eta \, f(\mu) =
 h(\mu) - \xi \eta \, f(\la) \,.
\ee
Hence we infer that $h(\la)$ is a linear function. Therefore,
the solution for $g(\la)$ looks as follows
\be\label{gsol2}
 g(\la)= \frac{\la}{\beta - \xi \eta \, \la} \,.
\ee
This function solves equation~\rf{G1} provided that
the following restrictions are imposed
\be\label{xibeta}
  \xi^2 = 1 \,, \qquad \beta = \eta - \xi/2 \,.
\ee

R--matrices corresponding to the cases 1), 2) and 3) are known
as the R--matrices of Yang, Baxter, and Zamolodchikovs, respectively.
Analysis presented above shows clearly that there are no other
solutions of the form~\rf{ans}. It is remarkable that the ansatz
\rf{ans} covers three out of the four known series of $sl_2$--invariant
regular R--matrices. It is therefore natural to study its
generalization that could be analysed in a similar way.

\section*{\S 3. On reduced Yang--Baxter equation}   %%%%%%%%%%%%%

Let us remind the main details of the approach developed in 
\cite{By} to analysis of $U_q(sl_2)$--invariant R--matrices
(we will take into account from the very beginning that $q=1$
in our case). Let the symbol $\lfloor t \rfloor$ denote the
entire part of a number~$t$. The subspace 
$W^{(s)}_n \subset V_s^{\otimes 3}$ for
\hbox{$n=0,1,\ldots,\lfloor 3s \rfloor$} is defined as a linear
span of highest weight vectors of spin $(3s - n)$, i.e.
\be\label{Wn}
  W^{(s)}_n = \{ \, \psi \in V_s^{\otimes 3} \quad \bigm| \quad
  S^+_{\one\two\three} \psi =0 \,, \quad
  S^z_{\one\two\three} \psi = (3s - n) \psi \, \} \,.
\ee
For a given R--matrix of form \rf{spec}, we construct a set of 
diagonal matrices $D^{(n)}(\la)$ as follows
\be\label{Dkk}
 D^{(n)}_{k k'}(\la) = \delta_{k k'} \, r_{2s-k}(\la) \,, \qquad
 \text{\rm where}\quad \biggl\{
 \begin{array}{ll}
  0 \leq k \leq n & \text{\rm for} \ \
    0 \leq n \leq 2s; \\ [0.5mm]
  n - 2s \leq k \leq 4s - n & \text{\rm for}\ \
      2s \leq n \leq \lfloor 3s \rfloor \,.
 \end{array}  \biggm.
\ee
We also introduce 
$\hat{D}^{(n)}(\la) \equiv A^{(s,n)} \ D^{(n)}(\la) \ A^{(s,n)}$,
where $A^{(s,n)}$ is a certain special matrix with properties
described below. Then the condition that the Yang--Baxter 
equation \rf{YB} is fulfilled on the subspace $W^{(s)}_n$ can 
be written as the following matrix equation
\begin{equation}\label{RYB}
 D^{(n)}(\la)  \ \hat{D}^{(n)}(\la+\mu) \ D^{(n)}(\mu)  =
 \hat{D}^{(n)}(\mu) \ D^{(n)}(\la+\mu) \ \hat{D}^{(n)}(\la)  \,,
\end{equation}
which we will call the {\em reduced} Yang--Baxter equation 
(of level~$n$). The initial equation (\ref{YB}) is equivalent 
to the system of reduced equations \rf{RYB} with 
\hbox{$n=0,1,\ldots,\lfloor 3s \rfloor$}.

The matrix $A^{(s,n)}$ which plays an important role in the
outlined approach has the following basic properties. Its 
entries are expressed (for $q=1$) in terms of 6--$j$ symbols
of the algebra $sl_2$ as follows (see also Appendix~A)
\be\label{A1}
 A^{(s,n)}_{k k^{\prime} }  =
 (-1)^{2s-n} \sqrt{(4s- 2k + 1)(4s - 2k' + 1)} \
 \FRW{s}{s}{s}{3s - n}{2s - k}{2s - k'} \,,
\ee
where $k,k'$ take values as in~\rf{Dkk}.
The matrix $A^{(s,n)}$ is orthogonal, symmetric, and coincides
with its own inverse ($t$ stands for matrix transposition)
\be\label{A3}
 A^{(s,n)} =  \bigl( A^{(s,n)} \bigr)^t =
  \bigl( A^{(s,n)} \bigr)^{-1} \,.
\ee

For the purpose of the present work we will need one more
property which we formulate as follows.
\begin{lem}\label{PM}
For all \hbox{$n=0,\ldots,{\lfloor 3s \rfloor}$} 
the following matrix relation holds
\be\label{AD0}
  A^{(s,n)} \, D^{(n)}_0 \, A^{(s,n)} =
    (-1)^n \, D^{(n)}_0 \, A^{(s,n)} \, D^{(n)}_0 \,,
\ee
where the diagonal matrix $D^{(n)}_0$ has the form
\be\label{D0}
 \bigl(D^{(n)}_0\bigr)_{k k'} = (-1)^k \delta_{k k'} \,,
\ee
and $k,k'$ take values as in~\rf{Dkk}.
\end{lem}
\proof
Let us write out matrix entries of \rf{AD0} taking into 
account that $A^{(s,n)}$ is symmetric:
\be\label{AD0b}
 \sum_m (-1)^m \, A^{(s,n)}_{k m} \, A^{(s,n)}_{k' m} =
 (-1)^{n+k+k'} A^{(s,n)}_{k k'} \,.
\ee
Now, taking into account formula \rf{A1}, it is easy 
to see that  relation  \rf{AD0b} can be reduced to
the Racah identity for 6--$j$ symbols (see, e.g.~\cite{BL})
\be\label{Racah}
 \sum_p (-1)^p  \, (2p+1) \,
 \FRW{r_1}{r_2}{r_3}{r_4}{l}{p}
 \FRW{r_1}{r_3}{r_2}{r_4}{l'}{p}
= (-1)^{l+l'} \, \FRW{r_3}{r_2}{r_1}{r_4}{l}{l'} \,,
\ee
where we have to set
$r_1=r_2=r_3=s$, $r_4=3s-n$, $l=2s-k$, $l'=2s-k'$, $p=2s-m$.
\qed\\[-0.5mm]

It is obvious from \rf{perm} that $D^{(n)}_0$ and
\hbox{$\hat{D}^{(n)}_0 \equiv A^{(s,n)} D^{(n)}_0 A^{(s,n)}$}
correspond to restriction of operators $\PP_{\one\two}$ and
$\PP_{\two\three}$ onto $W^{(s)}_n$. In particular, reduction
of  equation \rf{pp2} on the subspace $W^{(s)}_n$ leads
to the following relation
\begin{equation}\label{RYB0}
 D^{(n)}_0  \, A^{(s,n)} \, D^{(n)}_0 \, A^{(s,n)} \, D^{(n)}_0 =
 A^{(s,n)} \, D^{(n)}_0  \, A^{(s,n)} \, D^{(n)}_0 \,
 A^{(s,n)} \, D^{(n)}_0 \, A^{(s,n)} \,,
\end{equation}
correctness of which follows immediately from the statement of
Lemma~\ref{PM}. Another corollary of Lemma~\ref{PM} is that
$(-1)^n A^{(s,n)}$ corresponds to restriction of the operator
$\PP_{\one\three}=\PP_{\one\two}\PP_{\two\three}\PP_{\one\two}$
on~$W^{(s)}_n$.

\section*{\S 4. Partial analysis of a general ansatz}   %%%%%%%%%%%%%%%%
\subsection*{4.1 Derivation of equations}

Observe that any $sl_2$--invariant R--matrix of 
spin~$s \geq 1$ can be represented by the following ansatz
\be\label{ans2}
 R(\la) = \frac{1}{1+ f(\la)} \, \bigl( {\mathbb E} +
 f(\la) \, {\mathbb P} + g(\lambda) \, P^{2s-m} +
  \sum_{j=0}^{2s-m'}  \tilde{r}_j(\la) P^{j} \bigr) \,,
\ee
where $2 \leq m \leq 2s$ and $m < m'$ (if $m=2s$ then the last
sum in~\rf{ans2} is omitted). Below we will assume that
$g(\la)\neq 0$, since otherwise \rf{ans2} belongs to the known
case 1)~in~\S2. The regularity requirement imposes the condition
\be\label{reg2}
 f(0) = g(0) = \tilde{r}_j(0) = 0\,.
\ee

Let $\pi^{(m,n)}$ denote a matrix such that
$(\pi^{(m,n)})_{kk'} = \delta_{km} \delta_{k'm}$, $k=0,\ldots,n$. 
Then $\pi^{(m,n)}$ and
\hbox{$\hat{\pi}^{(m,n)} \equiv A^{(s,n)} \pi^{(m,n)} A^{(s,n)}$}
correspond to restriction of the operators
$P^{2s-m}_{\one\two}$ and $P^{2s-m}_{\two\three}$
on $W^{(s)}_n$. Notice that $\pi^{(m,n)}$ and $\hat{\pi}^{(m,n)}$
are projectors of rank~1.

For $n < m'$, the matrices $D^{(n)}(\la)$ and $\hat{D}^{(n)}(\la)$
corresponding to the R--matrix \rf{ans2} look as follows
\be\label{Dn}
\begin{aligned}
 D^{(n)}(\la) &= \frac{1}{1+ f(\la)} \,
     \bigl( \EE + f(\la) \, D^{(n)}_0 +
 \theta_{mn} \, g(\la) \, \pi^{(m,n)} \bigr) \,, \\
 \hat{D}^{(n)}(\la) &= \frac{1}{1+ f(\la)} \,
    \bigl( \EE + f(\la) \, \hat{D}^{(n)}_0 +
 \theta_{mn} \, g(\la) \, \hat{\pi}^{(m,n)} \bigr) \,,
\end{aligned}
\ee
where $\theta_{mn} = 0$ for $n < m$ and $\theta_{mn} = 1$ for 
$m \leq n < m'$.

The following observation is a key place of the present work:
analysis of the reduced Yang--Baxter equation \rf{RYB} for the 
ansatz \rf{Dn} is absolutely analogous (except for one special case) 
to analysis of  equation \rf{YB} for the ansatz \rf{ans}
given in~\S2. This observation is based on the following assertion:

\begin{lem}\label{PPrel2}
Relations \rf{pp1}--\rf{pp5} of Lemma~\ref{PPrel1} remain true
after the replacement
\be\label{zam1}
 \PP_{l} \to D^{(n)}_0 \,, \quad \PP_{l'} \to \hat{D}^{(n)}_0 \,,\quad
 P^0_{l} \to \pi^{(m,n)} \,, \quad P^0_{l'} \to \hat{\pi}^{(m,n)} \,,
\ee
as well as after the replacement
\be\label{zam2}
 \PP_{l} \to \hat{D}^{(n)}_0 \,, \quad \PP_{l'} \to D^{(n)}_0 \,,\quad
 P^0_{l} \to \hat{\pi}^{(m,n)} \,, \quad P^0_{l'} \to \pi^{(m,n)} \,.
\ee
The corresponding scalar constants $\xi$ and $\eta$ become 
dependent on $m$ and $n$:
\be\label{xieta2}
 \xi_m= (-1)^m \,, \qquad
 \eta_{m,n}=(-1)^n \, A^{(s,n)}_{m m}  \,.
\ee
\end{lem}
\proof
The analogues of relations \rf{pp1} follow from the definition of
matrices
$D^{(n)}_0$, $\hat{D}^{(n)}_0$, $\pi^{(m,n)}$, $\hat{\pi}^{(m,n)}$ 
and the property $(A^{(s,n)})^2 = \EE$. The analogue of relation
\rf{pp2} is identity \rf{RYB0} which we have established above. 
The analogues of relations \rf{pp3} can be reduced to the identity
\begin{equation}\label{pp3b}
 \pi^{(m,n)} \, A^{(s,n)} \, D^{(n)}_0  \, A^{(s,n)} \,
 D^{(n)}_0 \, A^{(s,n)}  = A^{(s,n)} \, D^{(n)}_0  \, A^{(s,n)} \,
 D^{(n)}_0 \,  A^{(s,n)} \, \pi^{(m,n)}  \,,
\end{equation}
which is easily verified with the help of  relation \rf{AD0}
and the analogues of relations~\rf{pp1}. The analogues of relations
\rf{pp4} and \rf{pp5} are derived as follows:
\begin{align*}
 \pi^{(m,n)} \, \hat{D}^{(n)}_0  \, \pi^{(m,n)} &=
 \pi^{(m,n)} \, A^{(s,n)}  D^{(n)}_0   A^{(s,n)} \, \pi^{(m,n)} \={AD0}
 (-1)^n \, \pi^{(m,n)} \, D^{(n)}_0  A^{(s,n)}  D^{(n)}_0 \, \pi^{(m,n)} \\
 &= (-1)^n \, \pi^{(m,n)} \, A^{(s,n)} \, \pi^{(m,n)} =
 (-1)^n \, A^{(s,n)}_{m m} \, \pi^{(m,n)} = \eta_{m,n} \, \pi^{(m,n)} \,, \\
\pi^{(m,n)} \, \hat{\pi}^{(m,n)}  \, \pi^{(m,n)} &=
 \pi^{(m,n)} \, A^{(s,n)} \, \pi^{(m,n)}  \, A^{(s,n)} \, \pi^{(m,n)}
 = \bigl(A^{(s,n)}_{m m}\bigr)^2 \, \pi^{(m,n)} =
 \eta_{m,n}^2 \, \pi^{(m,n)} \,, \\
\pi^{(m,n)} \, \hat{\pi}^{(m,n)}  \, D^{(n)}_0 &=
 \pi^{(m,n)} \, A^{(s,n)} \, \pi^{(m,n)}  \, A^{(s,n)}  D^{(n)}_0 \\
 &= A^{(s,n)}_{m m} \, \pi^{(m,n)}  \, A^{(s,n)}  D^{(n)}_0 \,
 \bigl( A^{(s,n)} \bigr)^2
 \={AD0} \eta_{m,n} \, \pi^{(m,n)} \, D^{(n)}_0
    A^{(s,n)}  D^{(n)}_0  A^{(s,n)} \\
 &= \xi_m \eta_{m,n} \, \pi^{(m,n)} \, A^{(s,n)}  D^{(n)}_0  A^{(s,n)}
 = \xi_m \eta_{m,n} \, \pi^{(m,n)} \, \hat{D}^{(n)}_0 \,.
\end{align*}
Let us emphasize that Lemma~\ref{PM} plays a key role in this proof.
\qed\\[-0.5mm]

The derivation of equations \rf{FGH} is based only on
relations of Lemma~\ref{PPrel1}. Therefore, Lemma~\ref{PPrel2} 
implies that for the R--matrix \rf{ans2} the reduced Yang--Baxter 
equation at levels $n<m'$ leads to the same equation~\rf{FGH}. 
The only difference is that the scalar coefficients 
(apart from $F_{\la,\mu}$) now depend on~$m$ and~$n$:
\be
\label{Fn}
 F_{\la,\mu} &=& f(\la) + f(\mu) - f(\la+\mu) \,, \\
\label{Gn}
 G^{(m,n)}_{\la,\mu} &=& \theta_{m,n} \, \Bigl(
 g(\la) + g(\mu) - g(\la+\mu) + \xi_{m} f(\la) g(\mu)
 + \xi_{m} \, g(\la) f(\mu) \\
\nn && +\ g(\la) g(\mu) + \eta_{m,n} g(\la) g(\mu) f(\la+\mu)
 + \eta^2_{m,n} g(\la) g(\mu) g(\la+\mu) \Bigr) \,, \\
\label{Hn}
 H^{(m,n)}_{\la,\mu} &=& \theta_{m,n} \, \Bigl(
 g(\la) f(\la+\mu) - f(\la) g(\la+\mu)
 + \xi_{m} \eta_{m,n} \, g(\la) f(\mu) g(\la+\mu) \Bigr) \,,
\ee
and matrices $\SF$, $\SG$, $\SH$, $\tilde{\SH}$ are given by 
formulae \rf{FGH0} after the substitution~\rf{zam1}, 
that is,
\begin{equation}\label{FGHR}
\begin{aligned}
\SF^{(m,n)} &= D^{(n)}_0 - \hat{D}^{(n)}_0 \,, \quad &
\SG^{(m,n)} &= \pi^{(m,n)} - \hat{\pi}^{(m,n)} \,, \\
\SH^{(m,n)} &= \pi^{(m,n)} \, \hat{D}^{(n)}_0 -
 D^{(n)}_0 \, \hat{\pi}^{(m,n)} \,, \quad &
\tilde{\SH}^{(m,n)} &=  \hat{D}^{(n)}_0 \, \pi^{(m,n)} -
 \hat{\pi}^{(m,n)} \, D^{(n)}_0  \,.
\end{aligned}
\end{equation}

\subsection*{4.2 Analysis of equations in the case $f(\la) = 0$}

Assuming that $g(\la)\neq 0$ in \rf{ans2}, let us consider first
the case $f(\la) = 0$.
In this case  equation \rf{FGH} at level $n=m$ is equivalent
to  equation $G_{\la,\mu}^{(m,m)} \SG^{(m,m)} =0$, i.e., 
to  equation \rf{G2} for~$g(\la)$, where $\eta$ has the form
\be\label{etann}
 \eta_{m,m} = (-1)^m \, A^{(s,m)}_{m m} =
 \frac{(2s)!}{(2s-m)!} \frac{(4s-2m+1)!}{(4s-m+1)!} \,.
%%  = \prod_{k=2s-m+1}^{4s-2m+1} \frac{k}{k+m}
\ee
For $2 \leq m \leq 2S$ and $S\geq 1$ we have
$|A^{(s,m)}_{m m}| < A^{(s,1)}_{1 1} = \frac{1}{2}$.
Therefore, $g(\la)$ is given by~\rf{gsol1}, where $\eta=\eta_{m,m}$.

Further analysis of the case $f(\la) = 0$ naturally leads  to a question:
is it possible for the ansatz \rf{ans2} to have $m'>(m+1)$? 
This is possible only if the already found function $g(\la)$ 
solves  equation \rf{G2} at level $n=m+1$, that is, only if
$\eta^2$ takes the same value for levels $n=m$ and $n=m+1$. 
According to \rf{xieta2}, the condition $\eta^2_{m,m}=\eta^2_{m,m+1}$ 
is equivalent to the requirement
\be\label{aann}
 | A^{(s,m)}_{m m} | = | A^{(s,m+1)}_{m m} | \,.
\ee
However, it is easy to derive from formula \rf{6j} that
\be\label{mm1}
  A^{(s,m+1)}_{m m} = \frac{m^2-m-3ms+s}{2s} \, A^{(s,m)}_{m m} \,.
\ee
Since $m^2-m-3ms+3s<0$ for $2 \leq m \leq 2s$, we infer that \rf{aann} 
cannot hold for these values of~$m$. Thus, we conclude that $m'=m+1$.

\subsection*{4.3 Analysis of equations in the case $f(\la) \neq 0$}

Let us now turn to the case $f(\la) \neq 0$. Equations \rf{FGH} at
levels $n=1,\ldots,m-1$ are equivalent to  equation
$F_{\la,\mu} \SF^{(m,n)} =0$, i.e., to  equation \rf{F1} for $f(\la)$. 
Therefore, without loss of generality we can choose $f(\la)=\la$.

In order to analyze  equation \rf{FGH} for $n\geq m$, it is
important to notice that the analogue of Lemma~\ref{LD} is
in general not true. That is, matrices \rf{FGHR} can be linearly
dependent. Observe that matrices $\SF^{(m,n)}$ and $\SG^{(m,n)}$ 
are symmetric and obviously linearly independent, whereas $\SH^{(m,n)}$
and $\tilde{\SH}^{(m,n)}$ are transposed to each other:
$\tilde{\SH}^{(m,n)} = \bigl(\SH^{(m,n)}\bigr)^t$. It turns out that
the following relations
\begin{align}
 \label{HHt}
 \tilde{\SH}^{(m,n)} &= \SH^{(m,n)} \,, \\
\label{HHbG}
 \SH^{(m,n)} + \tilde{\SH}^{(m,n)} &= \beta \, \SG^{(m,n)} \,,
\end{align}
where $\beta$ is a scalar constant, can hold only simultaneously. 
The case in which these relations do take place we will call
an {\em exceptional} one.
\begin{lem}\label{HHG}
For $m \geq 2$, each of  relations \rf{HHt} and \rf{HHbG}
holds only for $m=3$, $n=4$. In this case  relation \rf{HHbG}
holds in the following form
\be\label{HHGb}
 \SH^{(3,4)} + \tilde{\SH}^{(3,4)} = 2 \SG^{(3,4)} \,.
\ee
\end{lem}
The proof is given in Appendix~B.

In a generic case we have 
\hbox{$\SH^{(m,n)} \neq \tilde{\SH}^{(m,n)}$}.
Therefore, the antisymmetric matrix part of  equation \rf{FGH} 
imposes the condition
\be\label{H3}
 H^{(m,n)}_{\la,\mu} = H^{(m,n)}_{\mu,\la} \,,
\ee
and its symmetric part looks like following
\begin{equation}\label{FGH2}
 F_{\la,\mu} \, \SF^{(m,n)} +
 G^{(m,n)}_{\la,\mu} \, \SG^{(m,n)}  +
 \frac{1}{2}(H^{(m,n)}_{\la,\mu} + H^{(m,n)}_{\mu,\la}) \,
 ( \SH^{(m,n)}  + \tilde{\SH}^{(m,n)}) = 0 \,.
\end{equation}
If $\SF^{(m,n)}$, $\SG^{(m,n)}$ and
\hbox{$(\SH^{(m,n)} + \tilde{\SH}^{(m,n)})$} are linearly 
independent, then equations \rf{H3}--\rf{FGH2} lead to the 
system of functional equations~\rf{F1}--\rf{H1}. If, however,
$\SH^{(m,n)} + \tilde{\SH}^{(m,n)} =
 \beta \SG^{(m,n)} + \tilde{\beta} \SF^{(m,n)}$,
where $\tilde{\beta}\neq 0$, then \rf{FGH2} is equivalent
to the following system
\be \label{F3}
  2\, F_{\la,\mu} + \tilde{\beta} \,
  (H^{(m,n)}_{\la,\mu} + H^{(m,n)}_{\mu,\la}) &=& 0 \,, \\
\label{G3}
  2\, G^{(m,n)}_{\la,\mu} +
  \beta \, (H^{(m,n)}_{\la,\mu} + H^{(m,n)}_{\mu,\la}) &=& 0 \,.
\ee
Since our choice $f(\la)=\la$ has already ensured  equality
$F_{\la,\mu} =0$, we infer that  equations \rf{H3}, \rf{F3}--\rf{G3} 
lead again to system~\rf{F1}--\rf{H1}. Thus, we conclude
that analysis of a generic case is absolutely analogous to
analysis of the case 3) in~\S2.

Since, by Lemma~\ref{HHG}, the level \hbox{$n=m$} corresponds to
a generic case, the function $g(\la)$ is determined by
system~\rf{F1}--\rf{H1} uniquely and has the following form
\be\label{gsol3}
 g(\la)= \frac{\la}{\eta_{m,m} - \xi_m/2 - \xi_m \eta_{m,m} \, \la} \,,
\ee
where $\xi_{m}$ and $\eta_{m,m}$ are given by formulae~\rf{xieta2}.

Further analysis of the case $f(\la) = \la$ leads to a question:
is it possible for the ansatz \rf{ans2} to have \hbox{$m'>(m+1)$}?
If \hbox{$m \neq 3$}, then the level \hbox{$n=m+1$} corresponds to
a generic case. It is easy to check that the function \rf{gsol3} 
can satisfy system \rf{F1}--\rf{H1} for \hbox{$n=m+1$} only if
$\eta_{m,m}=\eta_{m,m+1}$. This is impossible, as it was shown in~\S4.2.
However, for $m=3$ this level corresponds to the exceptional case.
In this case  equation \rf{FGH2} is equivalent to the equation
\be\label{G4}
  G^{(3,4)}_{\la,\mu} + H^{(3,4)}_{\la,\mu} + H^{(3,4)}_{\mu,\la} = 0 \,,
\ee
whilst condition \rf{H3} is not imposed (because \rf{FGH} has no
antisymmetric part). Substituting \rf{gsol3} into \rf{G4}, it is
easy to verify that  equation \rf{G4} is true if 
\hbox{$(\eta_{3,3}-\eta_{3,4})(2\eta_{3,4}-1)=0$}. Interestingly, 
this condition is satisfied for all $s\geq \frac{3}{2}$, since,
according to \rf{etann} and \rf{mm1}, we have
\be\label{eta34}
 \eta_{3,4}= A^{(s,4)}_{3 3} = 1/2 \,.
\ee
Thus, for \hbox{$m = 3$} and for all $s\geq \frac{3}{2}$ we have
\hbox{$m'=4$} or \hbox{$m'=5$} in the ansatz~\rf{ans2}. Actually,
\hbox{$m'>5$} is possible only for $s=3$. Indeed, the level 
\hbox{$n=5$} corresponds to a generic case and hence a necessary
condition in order to have \hbox{$m'=6$} is the equality
\be\label{aann2}
 A^{(s,3)}_{3 3} = A^{(s,5)}_{3 3} \,.
\ee
However, it is not difficult to derive from \rf{6j} the following relation
\be\label{aann2b}
  A^{(s,5)}_{3 3} =
    \frac{10 s^2 - 32s + 21}{s(4s-7)} \, A^{(s,3)}_{3 3} \,,
\ee
which shows that \rf{aann2} can hold only if $(s-3)(6s-7)=0$. 
Finally, since $A^{(3,6)}_{3 3} \neq A^{(3,3)}_{3 3}$, 
we conclude that $m' \leq 6$ for~$s=3$.

\section*{\S5. Analysis of constant R--matrices}

We call $R \in {\rm End}\ V_{s}^{\otimes 2}$ a 
{\em constant} R--matrix if it solves the following Yang--Baxter equation
\be\label{YBc}
 R_{\one\two}\, R_{\two\three} \,
 R_{\one\two} =  R_{\two\three}\,
 R_{\one\two}\, R_{\two\three} \,.
\ee
We will consider $sl_2$--invariant R--matrices, i.e., those
that have the spectral decomposition,
\be\label{spec2}
 R = \sum_{j=0}^{2s}  r_j \, P^{j} \,,
\ee
where $r_j$ are scalar constants. In addition, we will assume that
\be\label{r2s}
 \qquad r_{2s} = 1\,.
\ee

The technique of analysing the spectral decomposition described
in \S\S3--4 is applicable to the case of constant R--matrices as 
well. In particular, the following remark explains why
condition \rf{r2s} is natural.
\begin{lem}\label{R2S}
There exist no nontrivial $sl_2$--invariant constant R--matrices
such that $r_{2s} = 0$.
\end{lem}
\proof
Indeed, if such an R--matrix exists, then by a suitable normalization
it can be brought to the form
\hbox{$R = P^{2s-m} + \sum_{j=0}^{2s-m-1}  r_j \, P^{j}$}, where
\hbox{$0<m \leq 2s$}. Then the corresponding reduced Yang--Baxter 
equation \rf{RYB} at the level \hbox{$n=m$} reads
$\bigl(A^{(s,m)} \, \pi^{m,m}\bigr)^3 =
    \bigl(\pi^{m,m} \, A^{(s,m)}\bigr)^3$,
which is equivalent to the relation
$(A^{(s,m)}_{m,m})^2 \, [ A^{(s,m)}, \pi^{m,m} ] =0$.
However, this relation cannot hold for $m>0$, because
$A^{(s,m)}_{k,m} \neq 0$ for all~$k$.
\qed\\[-0.5mm]

Observe that any $sl_2$--invariant constant R--matrix
of spin~$s \geq 1$ which satisfies~\rf{r2s} and has $r_{2s-1}\neq -1$
can be represented by the following ansatz
\be\label{ans5}
 R = \frac{1}{1+ f} \, \bigl( {\mathbb E} +
 f \, {\mathbb P} + g \, P^{2s-m} +
  \sum_{j=0}^{2s-m'}  \tilde{r}_j \, P^{j} \bigr) \,,
\ee
where $2 \leq m \leq 2s$ and $m < m'$ (if $m=2s$, then the last
sum in~\rf{ans5} is omitted).

Applying the same arguments as in \S3, we can use Lemma~\ref{PPrel2} 
to show that the reduced Yang--Baxter equation for the 
R--matrix \rf{ans5} at levels $n<m'$ is equivalent to the same 
equation \rf{FGH}, where the matrices $\SF$, $\SG$, $\SH$, $\tilde{\SH}$ 
are given by formulae \rf{FGHR}, and the scalar coefficients are 
obtained from \rf{Fn}--\rf{Hn} by converting the functions $f(\la)$ 
and $g(\la)$ into constants $f$ and $g$, i.e.:
\be
\label{Fn0}
 F &=& f \,, \\
\label{Gn0}
 G^{(m,n)} &=& \theta_{m,n} \, \bigl(
 g  + 2\xi_{m} f g +\ g^2 + \eta_{m,n} g^2 f
 + \eta^2_{m,n} g^3 \bigr) \,, \\
\label{Hn0}
 H^{(m,n)} &=& \theta_{m,n} \,
 \xi_{m} \eta_{m,n} \, g^2 f  \,.
\ee

Equation \rf{FGH} at levels $n < m$ is equivalent to  equation
$F \SF^{(m,n)} =0$, which can hold only for $f=0$. As a result, 
$G^{(m,n)}$ acquires the following form
\be
\label{Gn00}
 G^{(m,n)} &=& \theta_{m,n} \, \bigl(
 g  + g^2  + \eta^2_{m,n} g^3 \bigr)   \,,
\ee
and then  equation \rf{FGH} at levels $m\leq n < m'$ yields the 
equation $G^{(m,n)} \SG^{(m,n)} =0$, that is, the following quadratic 
equation on~$g$
\be
\label{gg}
 1  + g  + \eta^2_{m,n} g^2 = 0   \,.
\ee
Whence for $n=m$ we find $g=\frac{1}{2}(1 \pm \sqrt{1- 4 \eta^2_{m,m}})$,
where $\eta_{m,m}$ is given by formula~\rf{etann}.
Since, as it was shown in~\S4.2, we have $\eta^2_{m,m} \neq \eta^2_{m,m+1}$
for $2\leq m \leq 2s$, the obtained value of $g$ cannot
satisfy \rf{gg} for $n>m$. Thus, we conclude that $m'=m+1$ in~\rf{ans5}.

The ansatz \rf{ans5} covers not all $sl_2$--invariant constant
R--matrices of spin \hbox{$s \geq 1$} satisfying~\rf{r2s}. Namely,
if such an R--matrix has $r_{2s-1}= -1$, then it can be represented
by the following ansatz
\be\label{ans6}
 R =  {\mathbb P} + g \, P^{2s-m} +
  \sum_{j=0}^{2s-m'}  \tilde{r}_j \, P^{j}  \,,
\ee
where $2 \leq m \leq 2s$ and $m < m'$ (if $m=2s$, then the last
sum in~\rf{ans6} is omitted).

Using relations of  Lemma~\ref{PPrel2}, it is not difficult to
check that the reduced Yang--Baxter equation for the R--matrix \rf{ans6} 
at the level $n=m$ is equivalent to the following equation
\begin{equation}\label{GH}
 g^2 (1 + \eta_{m,m} \, g) \, \SG  +
 \xi_m \, g^2 \bigl(  \SH  + \tilde{\SH} \bigr) = 0 \,,
\end{equation}
where $\SG$, $\SH$, $\tilde{\SH}$ are given by formulae~\rf{FGHR}.
Since the level $n=m$ corresponds to a generic case (cf.~\S4.3), 
we infer that the only solution of \rf{GH} is $g=0$. That is, 
the ansatz \rf{ans6} is a solution for \rf{YBc} only if $g=r_j=0$.
Thus, we have shown that the permutation $\PP$ is a 
``rigid" solution, which does not admit a ``deformation" of its
spectral decomposition in the order $2s-2$ and lower orders.

\section*{Conclusion} %%%%%%%%%%%%%%%%%%%%%%%%%%%%%%%%%%%%%%%%%%%%

The results of analysis carried out in \S\S4--5 can be formulated
as the following restrictions on the structure of the spectral 
decomposition of R--matrices.

\begin{prop}\label{Hc0}
 Let $R$ be an $sl_2$--invariant solution of  equation \rf{YBc} 
on $V_s^{\otimes 3}$ for an integer or half--integer spin
$s \geq 1$, satisfying condition (\ref{r2s}).
Then either $r_{2s-1} = 1$ or \hbox{$r_{2s-1} = -1$}.

I. In the first case
\be\label{I0}
 R =   {\mathbb E}  + g \, P^{2s-m} +
  \sum_{j=0}^{2s-m-1}  \tilde{r}_j \, P^{j}  \,,
\ee
where $g$ is a solution of  equation~\rf{gg} and
$2 \leq m \leq 2s$. If $m<2s$, then $\tilde{r}_{2s-m-1}\neq 0$.

II. In the second case
\be\label{II0}
 R =  \PP  \,.
\ee

\end{prop}

\begin{prop}\label{Hc}
 Let $R(\la)$ be an $sl_2$--invariant solution of  equation \rf{YB} 
on $V_s^{\otimes 3}$ for an integer or half--integer spin
$s \geq 1$, satisfying conditions (\ref{run}).
Then either $r_{2s-1}(\la) = 1$ or
$r_{2s-1}(\la) = \frac{1- \gamma\la}{1+ \gamma\la}$.

I. In the first case
\be\label{I}
 R(\la) =   {\mathbb E}  + g(\la) \, P^{2s-m} +
  \sum_{j=0}^{2s-m-1}  \tilde{r}_j(\la) \, P^{j}  \,,
\ee
where $2 \leq m \leq 2s$. If $m<2s$, then
$\tilde{r}_{2s-m-1}(\la) \not\equiv 0$.
The function $g(\la)$ has the form
\be\label{Ig}
 g (\la) = b\, \frac{1-e^{\gamma\la}}%
  {e^{\gamma\la} - b^2} \,, \qquad
 b+b^{-1}= \frac{1}{\eta_{m,m}} \,,
\ee
where $\eta_{m,m}$ is given by~\rf{etann},
and $\gamma$ is some finite constant.

II. In the second case either
\be\label{IIa}
 R(\la) = \frac{1}{1+ \gamma\la} \, \bigl( {\mathbb E} +
 \gamma\la \, {\mathbb P} \bigr) \,,
\ee
or
\be\label{ans4}
 R(\la) = \frac{1}{1+ \gamma\la} \, \bigl( {\mathbb E} +
 \gamma\la \, {\mathbb P} +
 \frac{\la}{\eta_{m,m}(1 - (-1)^m \, \gamma\la) - \frac{(-1)^m}{2} } \,
 P^{2s-m} + \sum_{j=0}^{2s-m'}  \tilde{r}_j(\la) P^{j} \bigr) \,,
\ee
where $\eta_{m,m}$ is given by~\rf{etann},
$\gamma$ is some finite constant,
$2 \leq m \leq 2s$ and $m< m'$.
If $m<2s$, then $\tilde{r}_{2s-m'}(\la) \not\equiv 0$, and moreover
\begin{align*}
 {}&  m'=m+1 \,,\qquad\text{if\ } m\neq 3 \,, \\
 {}&  m' \leq 5 \,,\qquad\qquad  \text{if\ } m=3 \,, s\neq 3 \,, \\
 {}&  m' \leq 6 \,,\qquad\qquad  \text{if\ } m=3 \,, s = 3 \,.
\end{align*}
\end{prop}

The constant $\gamma$ can be set unity without loss of generality.

Let us make several brief remarks concerning  Propositions~1 and~2.

According to  Lemma~\ref{R2S}, if all coefficients in the 
spectral decomposition of $R(\la)$ tend to certain limit values
when $\la\to\infty$ in some direction in the complex plane, then
these values are finite. It follows from  Propositions~1 and~2
that the corresponding limit $R(\infty)$ has the form~\rf{I0}
only for solutions of the type~\rf{I}. In other cases 
we have $R(\infty)=\PP$.

What concerns  Proposition~2, we should remark
that for $s=3$ a solution with $m'=6$ really exists~\cite{Ke}:
\be\label{s3a}
 R(\lambda) = P^6 +
 \frac{1 {-} \lambda}{1 {+} \lambda } \, P^5 + P^4 +
  \frac{4 {-} \lambda}{4 {+} \lambda } \, P^3 + P^2 +
 \frac{1 {-} \lambda}{1 {+} \lambda } \, P^1 +
 \frac{1 {-} \lambda}{1 {+} \lambda } \,
 \frac{6 {-} \lambda}{6 {+} \lambda } \, P^0 \,.
\ee
It is easy to see that the coefficient of $P^3$ agrees with
formula~\rf{ans4}. Apart from this case, it is not known whether 
there exist R--matrices of the form \rf{ans4} with \hbox{$2<m<2s$}.

For $m=2$, the three highest order coefficients in~\rf{ans4},
\be\label{Rkrs}
 R(\lambda) = P^{2s} +
 \frac{1 {-} \lambda}{1 {+} \lambda } \, P^{2s-1} +
 \frac{1 {-} \lambda}{1 {+} \lambda } \,
 \frac{1 {-} \frac{2s}{2s-1}\lambda}{1 {+} \frac{2s}{2s-1}\lambda } \,
  P^{2s-2} + \ldots
\ee
coincide with the corresponding coefficients of the 
Kulish--Reshetikhin--Sklyanin R--matrix~\cite{KRS}. 
Let us mention that it follows from  Proposition~2 and
the results of~\cite{By} that only R--matrices of the form 
\rf{I} and \rf{Rkrs} can have $U_q(sl_2)$--invariant analogues.

%%%%%%%%%%%%%%%%%%%%%%%%%%%%%%%%%%%%%%%%%%%%%%
\appendix

%\section*{Appendix}
\subsection*{Appendix~A. Matrix $A^{(s,n)}$}
The expression \rf{A1} for entries of the matrix $A^{(s,n)}$ 
can be rewritten in a more explicit form:
\begin{align}
\label{6j}
 A^{(s,n)}_{k k'} &=
 F^s_k \, F^s_{k'} \, \sum_{l=6s-n-\min(k,k')}^{6s-\max(n,k+k')}
 (-1)^l (l{+}1)!
 \, \Bigl(  (l{-}4s{+}k)! \, (l{-}4s{+}k^\prime)! \Bigm. \\
\nonumber
 & \times   \Bigm.
 (l{-}6s{+}n{+}k)! \,(l{-}6s{+}n{+}k^\prime)! \,
 (6s{-}n{-}l)! \, (6s{-}k{-}k^\prime{-}l)! \,
 (8s{-}n{-}k{-}k^\prime{-}l)! \Bigr)^{-1} \,,
\end{align}
where $k,k'$ take values as in~\rf{Dkk} and
\be
\label{Fsk}
 F^s_{k} = (2s {-} k)! \, \Bigl( \frac{(k)!\,
 (n{-}k)! \, (2s{-}n{+}k)! \, (4s{-}n{-}k)!}%
 {(4s{-}k{+}1)! \, (6s{-}n{-}k{+}1)!} \Bigr)^{\frac{1}{2}} \,.
\ee
The summation in (\ref{6j}) is taken over those $l$ for which
the arguments of factorials are nonnegative and it is 
understood that $0!=1$.

\subsection*{Appendix~B. Proof of Lemma~\ref{LD} and Lemma~\ref{HHG}}

{\em Proof} of Lemma~\ref{LD}:
Using the relations ($a,b=1,2,3$)
\be\label{tr}
 \tr_a \EE_{a} = 2s+1 \,, \qquad
 \tr_a P^j_{ab} = \frac{2j+1}{2s+1} \, \EE_{b} \,, \qquad
 \tr_a \PP_{ab} = \EE_{b} \,,
\ee
we take the trace over the third tensor component of $\SF$, $\SG$, 
$\SH$ and $\tilde{\SH}$, which yields:
\be
 \tr_3 \SF = \eta^{-1} \PP - \EE \,, \qquad
 \tr_3 \SG = \eta^{-1} P^0 - \eta \EE \,, \qquad
 \tr_3 \SH = \tr_3 \tilde{\SH} = P^0 - \eta \PP \,.
\ee
Since $\EE$, $\PP$, $P^0$ are linearly independent for $s\geq 1$, 
we conclude that $\SF$, $\SG$, $\SH$ and $\tilde{\SH}$ can be 
linearly dependent only if the following equality holds
\be\label{lindep1}
 \eta^2 \, \SF - \eta \, \SG + \alpha \, \SH +
  \tilde{\alpha} \, \tilde{\SH} = 0 \,, 
  \qquad \alpha + \tilde{\alpha} =1 \,.
\ee
Multiply \rf{lindep1} by $P^0_{\one\two}$ from the left, taking into 
account  relations~\rf{pp1}, and take the trace over the first 
tensor component. Using again linear independence of $\EE$, $\PP$, 
$P^0$, we infer that \rf{lindep1} can hold only if $\alpha=\xi \eta$. 
Multiplying \rf{lindep1} by $P^0_{\one\two}$ from the right, we infer
analogously that $\tilde{\alpha} = \xi \eta$. Thus, \rf{lindep1} can 
hold only if $\xi \eta = 1/2$, which is impossible as seen from~\rf{xieta}.
\qed\\[-0.5mm]

\noindent
{\em Proof} of Lemma~\ref{HHG}:
Let us write out entries of the matrices $\SH$, $\tilde{\SH}$ 
and $\SG$ explicitly:
\begin{align}
\label{Hkk}
  H _{kk'} &= (-1)^{m+k'} \, \delta_{km} \, A^{(s,n)}_{kk'}
 - (-1)^k \, A^{(s,n)}_{km} \, A^{(s,n)}_{mk'} \,, \\
 \tilde{H}_{kk'} &=
   (-1)^{m+k} \, \delta_{km} \, A^{(s,n)}_{kk'}
  - (-1)^{k'} \, A^{(s,n)}_{km} \, A^{(s,n)}_{mk'} \,, \\
\label{Gkk}
 G_{kk'} &= \delta_{km} \, \delta_{k' m} \,  -
    A^{(s,n)}_{km} \, A^{(s,n)}_{mk'} \,.
\end{align}
Recall that $m\geq 2$ and $k,k'=0,1,\ldots,n$. Comparing
\rf{Hkk}--\rf{Gkk} for $k=k'=0$, we notice that \rf{HHbG} 
can hold only for $\beta=2$. Further, considering
\rf{Hkk}--\rf{Gkk} for $k,k'\neq m$, it is easy to see that
each of  relations \rf{HHt} and \rf{HHbG} can hold only if
\be\label{AmAm}
   A^{(s,n)}_{km} \, A^{(s,n)}_{mk'} = 0 \,,
\ee
for all values of $k, k'$ such that $k,k'\neq m$ and
$(-1)^k + (-1)^{k'} \neq 2$. In particular, \rf{AmAm} must
hold for $k=k'=1$, which implies that $A^{(s,n)}_{1m}=0$.
As can be seen from \rf{6j}, the latter equality is possible only
if the following condition is satisfied (for $m<n$)
\be\label{mns1}
 2m^2 - 2m + n^2 - n = 8ms - 6ns \,.
\ee
Observe that $A^{(s,n)}_{kn} \neq 0$ for all~$k$.
Therefore, \rf{AmAm} implies that $m \neq n$, and also
that $n$ is an even number. Furthermore, if $m \neq n-1$, then
\rf{AmAm} must hold for $k=k'=n-1$, which implies that
$A^{(s,n)}_{1,n-1}=0$. But, as can be inferred again from \rf{6j},
this is possible only if the following condition is fulfilled
\be\label{mns2}
 m^2 - m = 4ms - ns \,.
\ee
It is easy to see that the conditions \rf{mns1} and \rf{mns2} 
are incompatible because they imply the equality $n^2 - n + 4ns =0$.
Thus, the only remaining possibility is the case in which $m=n-1$. 
In this case condition \rf{mns1} is satisfied only for $n=4$. 
A direct check shows that  relations \rf{HHt} and \rf{HHGb}
indeed hold for $m=3$, $n=4$.
\qed

%%%%%%
\newcommand{\my}[7]{{#1:} {\em #2.---} {#3} {\bf #4} {(#5),} {#6}}
\small  \setlength{\itemsep}{-3pt}

\end{document}